\documentclass[12pt]{article}

\usepackage{graphicx,psfrag}
\usepackage{psfrag,eepic,epsfig}
\usepackage{sectsty}
\usepackage{setspace}
\usepackage{color}

\usepackage{setspace}
\usepackage{amsmath, epsfig, cite, ifpdf}
\usepackage{amssymb}
\usepackage{amsfonts}
\usepackage{latexsym}
\usepackage{graphicx,psfrag,eepic,rawfonts}
\usepackage{enumerate}
\usepackage{indentfirst}
\input{epsf}

\newtheorem{thm}{Theorem}[section]

\newtheorem{cor}[thm]{Corollary}

\newtheorem{lem}[thm]{Lemma}

\newtheorem{problem}[thm]{Problem}
\newtheorem{exam}[thm]{Example}
\newtheorem{remark}[thm]{Remark}

\newcommand{\qed}{{\hfill\rule{4pt}{7pt}}}
\def\pf{\noindent {\it Proof.} }

\numberwithin{equation}{section}

\makeatletter \@addtoreset{equation}{section} \makeatother

\title {\bf Skew-spectra and skew energy of various products of
graphs\footnote{Supported by NSFC and the ``973" program. }}

\author{
{\small  Xueliang Li, Huishu Lian}\\
{\small Center for Combinatorics and LPMC-TJKLC}\\
{\small Nankai University, Tianjin 300071, P.R. China}\\
{\small E-mail: lxl@nankai.edu.cn; lhs6803@126.com}
   }
\date{}
\begin{document}

\maketitle

\begin{abstract}
Given a graph $G$, let $G^\sigma$ be an oriented graph of $G$ with
the orientation $\sigma$ and skew-adjacency matrix $S(G^\sigma)$.
Then the spectrum of $S(G^\sigma)$ consisting of all the eigenvalues of
$S(G^\sigma)$ is called the skew-spectrum of $G^\sigma$, denoted by
$Sp(G^\sigma)$. The skew energy of the oriented graph $G^\sigma$,
denoted by $\mathcal{E}_S(G^\sigma)$, is defined as the sum of the
norms of all the eigenvalues of $S(G^\sigma)$. In this paper,
we give orientations of the Kronecker product $H\otimes G$ and the strong
product $H\ast G$ of $H$ and $G$ where $H$ is a bipartite graph and $G$
is an arbitrary graph. Then we determine the skew-spectra of the resultant
oriented graphs. As applications, we construct new families of oriented
graphs with maximum skew energy. Moreover, we consider the skew energy of
the orientation of the lexicographic product $H[G]$ of a bipartite graph $H$
and a graph $G$.


\end{abstract}

\section{Introduction}

Let $G$ be a simple undirected graph with vertex set
$V(G)=\{v_1,v_2,\ldots,v_n\}$, and let $G^\sigma$ be an oriented
graph of $G$ with the orientation $\sigma$, which assigns to each
edge of $G$ a direction so that the induced graph $G^\sigma$ becomes
an oriented graph or a directed graph. Then $G$ is called the underlying
graph of $G^\sigma$. The {\it skew-adjacency matrix} of $G^\sigma$ is
the $n\times n$ matrix $S(G^\sigma)=[s_{ij}]$, where $s_{ij}=1$ and
$s_{ji}=-1$ if $\langle v_i,v_j\rangle$ is an arc of $G^\sigma$,
otherwise $s_{ij}=s_{ji}=0$.  It is easy to see that $S(G^\sigma)$ is a
skew-symmetric matrix, and thus all its eigenvalues are purely imaginary
numbers or 0, which form the spectrum of $S(G^\sigma)$ and are said to be the
{\it skew-spectrum} $Sp(G^\sigma)$ of $G^\sigma$.

The concept of the energy of a simple undirected graph was introduced by
Gutman in \cite{G}. Then Adiga, Balakrishnan and So in \cite{ABC} generalized
the energy of an undirected graph to the skew energy of an oriented graph.
Formally, the {\it skew energy} of an oriented graph $G^\sigma$ is defined
as the sum of the absolute values of all the eigenvalues of $S(G^\sigma)$,
denoted by $\mathcal{E}_S(G^\sigma)$. Most of the results on the skew energy
are collected in our recent survey \cite{LL}, among which the problem about
the maximum skew energy has been paid more attention.

In \cite{ABC}, Adiga, Balakrishnan and So derived that for any oriented
graph $G^\sigma$ with order $n$ and maximum degree $\Delta$,
$\mathcal{E}_S(G^\sigma)\leq n\sqrt{\Delta}$. They also showed that the
equality holds if and only if $S(G^\sigma)^{T}S(G^\sigma)=\Delta I_n$,
which implies that $G^\sigma$ is $\Delta$-regular. Among all oriented
graphs with order $n$ and maximum degree $\Delta$, the skew
energy $n\sqrt{\Delta}$ is called the maximum skew energy. Naturally,
they proposed the following problem:
\begin{problem}\label{prob1}
Which $k$-regular graphs on $n$ vertices have orientations $G^\sigma$ with
$\mathcal{E}_S(G^\sigma)= n\sqrt{k}$, or equivalently, $S(G^\sigma)^{T}S(G^\sigma)=k I_n$ ?
\end{problem}

For $k=1,2,3,4$, all $k$-regular graphs which have orientations $G^\sigma$
with $\mathcal{E}_S(G^\sigma)= n\sqrt{k}$ were characterized, see
\cite{ABC,GX,CLL}. Other families of oriented regular graphs with the maximum
skew energy were also obtained. Tian in \cite{Tian} gave the orientation of
the hypercube $Q_k$ such that the resultant oriented graph has maximum
skew energy. In \cite{ABC}, a family of oriented graphs with maximum skew
energy was constructed by considering the Kronecker product of graphs. To be
specific, let $G^{\sigma_1}_1$, $G^{\sigma_2}_2$, $G^{\sigma_3}_3$ be the oriented
graphs of order $n_1$, $n_2$, $n_3$ with skew-adjacency matrices $S_1$, $S_2$,
$S_3$, respectively. Then the Kronecker product matrix $S_1\otimes S_2\otimes S_3$
is also skew-symmetric and is in fact the skew-adjacency matrix of an oriented graph
of the Kronecker product $G_1\otimes G_2\otimes G_3$. Denote the corresponding
oriented graph by $G^{\sigma_1}_1\otimes G^{\sigma_2}_2\otimes G^{\sigma_3}_3$.
The following result was obtained.
\begin{thm}\cite{ABC}\label{odd-Kron}
Let $G^{\sigma_1}_1$, $G^{\sigma_2}_2$, $G^{\sigma_3}_3$ be the oriented regular
graphs of order $n_1$, $n_2$, $n_3$ with maximum skew energies $n_1\sqrt{k_1}$,
$n_2\sqrt{k_2}$, $n_3\sqrt{k_3}$, respectively. Denote by $S_1$, $S_2$ and $S_3$ the
skew-adjacency matrices of $G^{\sigma_1}_1$, $G^{\sigma_2}_2$ and $G^{\sigma_3}_3$,
respectively. Then the oriented graph $G^{\sigma_1}_1\otimes G^{\sigma_2}_2\otimes G^{\sigma_3}_3$
has maximum skew energy $n_1n_2n_3\sqrt{k_1k_2k_3}$.
\end{thm}

It should be noted that the above Kronecker product of oriented graphs is
naturally defined, but the product requires $3$ or an odd number of oriented graphs.

Moreover, Cui and Hou in \cite{CH} gave an orientation $(P_m\square G)^o$ of the
Cartesian product $P_m\square G$, where $P_m$ is a path of order $m$ and $G$ is an
arbitrary graph. They computed the skew-spectra of $(P_m\square G)^o$, and by applying
this result they constructed a family of oriented graphs with maximum skew energy.
Then we in \cite{CLL3} extended their results to the oriented graph $(H\square G)^o$
where $H$ is an arbitrary bipartite graph, and thus a larger family of oriented graphs 
with maximum skew energy was obtained.
\begin{thm} \cite{CLL3}\label{Cartesian}
Let $H^\tau$ be an oriented $\ell$-regular bipartite graph on $m$ vertices with maximum
skew energy $\mathcal{E_S}(H^\tau)=m\sqrt{\ell}$ and $G^\sigma$ be an oriented
$k$-regular graph on $n$ vertices with maximum skew energy $\mathcal{E_S}(G^\sigma)=n\sqrt{k}$.
Then the oriented graph $(H^\tau\Box G^\sigma)^o$ of $H\Box G$ has the maximum skew
energy $\mathcal{E_S}((H^\tau\square G^\sigma)^o)=mn \sqrt{\ell+k}$.
\end{thm}

In this paper, we consider other products of graphs, including the Kronecker product
$H\otimes G$, the strong product $H\ast G$ and the lexicographic product $H[G]$, where
$H$ is a bipartite graph and $G$ is an arbitrary graph. In Subsection \ref{otimes},
We first give an orientation of $H\otimes G$, and then determine the skew-spectra of
the resultant oriented graph. As an application, we construct a new family of oriented
graphs with maximum skew energy. Subsection \ref{ast} is used to orient the graph
$H\ast G$, determine the skew-spectra of the resultant oriented graph and construct
another new family of oriented graphs with maximum skew energy. Finally we consider
the skew energy of the orientation of the lexicographic product $H[G]$ of $H$ and $G$
in Subsection \ref{lexi}.

In the sequel of this paper, it will be seen that there is no limitation of the number
of oriented graphs in our Kronecker product, and the oriented graphs that we will construct
have smaller order than the previous results under the same regularity.

\section{Main results}

We first recall some definitions. Let $H$ be a graph of order $m$ and $G$ be a graph
of order $n$. The Cartesian product $H\square G$ of $H$ and $G$ has vertex set
$V(H)\times V(G)$, where $(u_1,v_1)$ is adjacent to $(u_2, v_2)$ if and only if
$u_1=u_2$ and $v_1$ is adjacent to $v_2$ in $G$, or $u_1$ is adjacent to $u_2$ in $H$
and $v_1=v_2$. The Kronecker product $H\otimes G$ of $H$ and $G$ is a graph with vertex
set $V(H)\times V(G)$ and where $(u_1,v_1)$ and $(u_2, v_2)$ are adjacent if $u_1$ is
adjacent to $u_2$ in $H$ and $v_1$ is adjacent to $v_2$ in $G$. The strong product
$H\ast G$ of $H$ and $G$ is a graph with vertex set $V(H)\times V(G)$; two distinct
pairs $(u_1,v_1)$ and $(u_2, v_2)$ are adjacent in $H\ast G$ if $u_1$ is equal or
adjacent to $u_2$, and $v_1$ is equal or adjacent to $v_2$. The lexicographic product
$H[G]$ of $H$ and $G$ has vertex set $V(H)\times V(G)$ where $(u_1,v_1)$ is adjacent
to $(u_2, v_2)$ if and only if $u_1$ is adjacent to $u_2$ in $H$, or $u_1=u_2$ and
$v_1$ is adjacent to $v_2$ in $G$.

It can be verified that the Cartesian product $H\square G$, the Kronecker product
$H\otimes G$, the strong product $H\ast G$ are commutative, that is, $H\square G=G\square H$,
$H\otimes G=G\otimes H$ and $H\ast G=G\ast H$. But the lexicographic product $H[G]$ may
not be the same as $G[H]$. Moreover, the two graphs $H\square G$ and $H\otimes G$ are
edge-disjoint and $E(H\ast G)=E(H\square G)\cup E(H\otimes G)$. Finally, we point out
that if $H$ is a bipartite graph, then $H\otimes G$ is also bipartite.

In what follows, we always assume that $H$ is a bipartite graph on $m$ vertices
with bipartite $(X,Y)$ where $|X|=m_1$ and $|Y|=m_2$ and $G$ is a graph on $n$ vertices.
Let $H^\tau$ be an arbitrary oriented graph of $H$ and $G^\sigma$ be an arbitrary oriented
graph of $G$. Let $S_1$ and $S_2$ be the skew-adjacency matrices of $H^\tau$ and $G^\sigma$,
respectively. Giving the labeling of the vertices of $H$ such that the vertices of $X$ are
labeled first. Then the skew-adjacency matrix $S_1$ can be formulated as
$\left(\begin{array}{cc} 0& A\\ -A^T& 0\end{array}\right)$, where $A$ is an $m_1\times m_2$
matrix and $m_1+m_2=m$. Let $S'_1=\left(\begin{array}{cc} 0& A\\ A^T& 0\end{array}\right)$.
Note that $S_1$ is skew-symmetric and $S'_1$ is symmetric. It is easy to see that
$S_1S_1^T=S'_1(S'_1)^T$, and thus $S_1$ and $S_2$ have the same singular values.

\subsection{The orientation of $H\otimes G$}\label{otimes}

We first give an orientation of $H\otimes G$. For any two adjacent vertices $(u_1,v_1)$
and $(u_2,v_2)$, $u_1$ and $u_2$ must be in different parts of the bipartition of vertices of $H$ and assume
that $u_1\in X$. Then there is an arc from $(u_1,v_1)$ to $(u_2,v_2)$ if $\langle u_1,u_2\rangle$
is an arc of $H$ and $\langle v_1,v_2\rangle$ is an arc of $G$, or $\langle u_2,u_1\rangle$
is an arc of $H$ and $\langle v_2,v_1\rangle$ is an arc of $G$; otherwise there is an arc
from $(u_2,v_2)$ to $(u_1,v_1)$. Denote by $(H^\tau \otimes G^\sigma)^o$ the resultant
oriented graph and by $S$ its skew-adjacency matrix. For the skew-spectrum of
$(H^\tau \otimes G^\sigma)^o$, we obtain the following result.
\begin{thm}\label{spectrum}
Let $H^\tau$ be an oriented bipartite graph of order $m$ and let the
skew-eigenvalues of $H^\tau$ be the non-zero values $\pm \mu_1i$, $\pm
\mu_2i$, $\dots$, $\pm \mu_ti$ and $m-2t$ $0$'s. Let $G^\sigma$ be an oriented
graph of order $n$ and let the skew-eigenvalues of $G^\sigma$ be the non-zero
values $\pm \lambda_1i$, $\pm\lambda_2i$, $\dots$, $\pm \lambda_ri$ and $n-2r$ $0$'s.
Then the skew-eigenvalues of the oriented graph $(H^\tau\otimes G^\sigma)^o$ are
$\pm \mu_j\lambda_ki$ with multiplicities $2$, $j=1,\dots,t$, $k=1,\dots,r$,
and $0$ with multiplicities $mn-4rt$.
\end{thm}
\pf With suitable labeling of the vertices of $H\otimes G$, the skew-adjacency matrix
$S$ of $(H^\tau \otimes G^\sigma)^o$ can be formulated as $S=S'_1\otimes S_2$. We first
compute the singular values of $S$. Note that $S^T=-S'_1\otimes S_2$. Then
$$SS^T= (S'_1\otimes S_2)(-S'_1\otimes S_2)=-(S'_1)^2\otimes S^2_2.$$
It follows that the eigenvalues of $SS^T$ are $\mu(S'_1)^2\cdot \lambda(S_2)^2$, where
$\mu(S'_1)$ is the eigenvalues of $S'_1$ and $\lambda(S_2)i$ is the eigenvalues of $S_2$.
That is to say, the eigenvalues of $SS^T$ are $\mu(H^\tau)^2\cdot \lambda(G^\sigma)^2$,
where $\mu(H^\tau)i \in Sp(H^\tau)$ and $\lambda(G^\sigma)i\in Sp(G^\sigma)$.
From this, it immediately follows what we want. The proof is now complete. \qed

The above theorem can be used to yield a family of oriented graphs with maximum skew energy.
The following lemma was obtained in \cite{ABC}.
\begin{lem}\cite{ABC}\label{max}
Let $G^\sigma$ be an oriented graph of $G$ with order $n$ and maximum degree $\Delta$.
Then $\mathcal{E}_S(G^\sigma)\leq n\sqrt{\Delta}$, where the equality holds if and only if $S(G^\sigma)^{T}S(G^\sigma)=\Delta I_n$.
\end{lem}

\begin{thm}\label{max-fami1}
Let $H^\tau$ be an oriented $k$-regular bipartite graph of order $m$ with maximum
skew energy $m\sqrt{k}$. Let $G^\sigma$ be an oriented $l$-regular graph of order $n$
and the maximum skew energy $n\sqrt{l}$. Then $(H^\tau\otimes G^\sigma)^o$ is an oriented $kl$-regular
bipartite graph and has the maximum skew energy $\mathcal{E}_S((H^\tau\otimes G^\sigma)^o)=mn\sqrt{kl}$.
\end{thm}
\pf By the definition of the Kronecker product, it is easy to find that
$H\otimes G$ is a $kl$-regular bipartite graph with $mn$ vertices.
Let $S(H^\tau)=\left(\begin{array}{cc} 0& A\\ -A^T& 0\end{array}\right)$ be the
skew-adjacency matrix of $H^\tau$ and $S(G^\sigma)$ be the skew-adjacency matrix of
$G^\sigma$. Then by Lemma \ref{max}, we have $S(H^\tau)^TS(H^\tau)=kI_m$
and $S(G^\sigma)^TS(G^\sigma)=lI_n$. From Theorem \ref{spectrum}, the skew-adjacency matrix
$S$ of $(H^\tau\otimes G^\sigma)^o$ can be written as $S=S'(H^\tau)\otimes S(G^\sigma)$, where
$S'(H^\tau)=\left(\begin{array}{cc} 0& A\\ A^T& 0\end{array}\right)$. Note that
$S'(H^\tau)^TS'(H^\tau)=S(H^\tau)^TS(H^\tau)=kI_m$. It follows that
\begin{equation*}
\begin{split}
S^TS&=\left(S'(H^\tau)\otimes S(G^\sigma)\right)^T\left(S'(H^\tau)\otimes S(G^\sigma)\right)\\
&=\left((S'(H^\tau))^TS'(H^\tau)\right)\otimes \left(S(G^\sigma))^T(S(G^\sigma)\right)\\
&=kI_m\otimes lI_n=klI_{mn}.
\end{split}
\end{equation*}
By Lemma \ref{max}, the oriented graph $(H^\tau\otimes G^\sigma)^o$ has the maximum
skew energy $\mathcal{E}_S((H^\tau\otimes G^\sigma)^o)=mn\sqrt{kl}$. We thus complete
the proof of this theorem. \qed

Let $H^\tau$ be an oriented bipartite graph with maximum skew energy. Let $G_1^{\sigma_1}$
and $G_2^{\sigma_2}$ be any two oriented graphs with the maximum skew energies. By the above
theorem, the oriented graph $(H^\tau\otimes G_1^{\sigma_1})^o$ is bipartite and has the
maximum skew energy. Therefore, the Kronecker product $H\otimes G_1\otimes G_2$ can be
oriented as $((H^\tau\otimes G_1^{\sigma_1})^o\otimes G_2^{\sigma_2})^o$, abbreviated as
$(H^\tau\otimes G_1^{\sigma_1}\otimes G_2^{\sigma_2})^o$, which is also bipartite and has
the maximum skew energy. The process is valid for any positive integral number of
oriented graphs. Then the following corollary is immediately implied.
\begin{cor}\label{Kron}
Let $H^\tau$ be an oriented $k$-regular bipartite graph of order $m$ with maximum skew
energy $m\sqrt{k}$. Let $G_i^{\sigma_i}$ be an oriented $l_i$-regular graph of order $n_i$
with maximum skew energy $n_i\sqrt{l_i}$ for $i=1,2,\dots,s$ and any positive integer $s$.
Then the oriented graph $(H^\tau\otimes G_1^{\sigma_1}\otimes \cdots\otimes G_s^{\sigma_s})^o$
has the maximum skew energy $mn_1n_2\cdots n_s\sqrt{k\,l_1l_2\cdots l_s}$.
\end{cor}

\begin{remark}
In Corollary \ref{Kron}, the value $s$ can be any positive integer. If $s=2$, then the
oriented graph $(H^\tau\otimes G_1^{\sigma_1}\otimes G_2^{\sigma_2})^o$ has the maximum
skew energy $mn_1n_2\sqrt{k\,l_1l_2}$. Recall the orientation in Theorem \ref{odd-Kron},
which illustrates that the oriented graph $H^\tau\otimes G_1^{\sigma_1}\otimes G_2^{\sigma_2}$
also has the maximum skew energy $mn_1n_2\sqrt{k\,l_1l_2}$. In fact, this two orientations are identical.

Let $S_0=\left(\begin{array}{cc} 0& A\\ -A^T& 0\end{array}\right)$ be the skew-adjacency
matrix of $H^\tau$ and $S'_0=\left(\begin{array}{cc} 0& A\\ A^T& 0\end{array}\right)$.
Let $S_1$ and $S_2$ be the skew-adjacency matrices of $G_1^{\sigma_1}$ and $G_2^{\sigma_2}$.
Then the oriented graph $H^\tau\otimes G_1^{\sigma_1}\otimes G_2^{\sigma_2}$ has the
skew-adjacency matrix $\left(\begin{array}{cc} 0& A\otimes S_1\otimes S_2\\ -A^T\otimes S_1\otimes S_2 & 0\end{array}\right)$. The oriented graph $(H^\tau\otimes G_1^{\sigma_1})^o$ has the skew-adjacency matrix
$$S'_1\otimes S_2=\left(\begin{array}{cc} 0& A\otimes S_1\\ A^T\otimes S_1& 0\end{array}\right).$$
It follows that the skew-adjacency matrix of $(H^\tau\otimes G_1^{\sigma_1}\otimes G_2^{\sigma_2})^o$ is
$$\left(\begin{array}{cc} 0& A\otimes S_1\\ -A^T\otimes S_1& 0\end{array}\right)\otimes S_2=\left(\begin{array}{cc} 0& A\otimes S_1\otimes S_2\\ -A^T\otimes S_1\otimes S_2 & 0\end{array}\right),$$
which is the same as that of $H^\tau\otimes G_1^{\sigma_1}\otimes G_2^{\sigma_2}$.

In fact, for any even $s$, the oriented graph obtained in Corollary \ref{Kron} is identical
to the one obtained in Theorem \ref{odd-Kron}.
\end{remark}

\subsection{The orientation of $H\ast G$}\label{ast}

Now we consider the strong product $H\ast G$ of a bipartite graph $H$ and a graph $G$,
Let $H^\tau$ be an oriented graph of $H$ and $G^\sigma$ be an oriented graph of $G$.
Since the edge set of $H\ast G$ is the disjoint-union of the edge sets of $H\square G$
and $H\otimes G$, there is a natural orientation of $H\ast G$ if $H\square G$ and
$H\otimes G$ have been given orientations.

First recall the orientation $(H^\tau\square G^\sigma)^o$ of the Cartesian product
$H\square G$ given in \cite{CLL3}. For any two adjacent matrices $(u_1,v_1)$ and
$(u_2,v_2)$, we give it an orientation as follows. When $u_1=u_2\in X$, there is an
arc from $(u_1,v_1)$ to $(u_2,v_2)$ if $\langle v_1,v_2\rangle$ is an arc of $G^\sigma$
and an arc from $(u_2,v_2)$ to $(u_1,v_1)$ otherwise. When $u_1=u_2\in Y$, there is
an arc from $(u_1,v_1)$ to $(u_2,v_2)$ if $\langle v_2,v_1\rangle$ is an arc of $G^\sigma$
and an arc from $(u_2,v_2)$ to $(u_1,v_1)$ otherwise. When $v_1=v_2$, there is an arc
from $(u_1,v_1)$ to $(u_2,v_2)$ if $\langle u_1,u_2\rangle$ is an arc of $H^\tau$ and
an arc from $(u_2,v_2)$ to $(u_1,v_1)$ otherwise. Let $\overline{S}$ be the skew-adjacency
matrix of $(H^\tau\square G^\sigma)^o$.

Now we give an orientation of $H\ast G$ such that the arc set of the resultant oriented
graph is the disjoint-union of the arc sets of $(H^\tau\square G^\sigma)^o$ and
$(H^\tau\otimes G^\sigma)^o$. Denote by $(H^\tau\ast G^\sigma)^o$ this resultant
oriented graph and by $\widehat{S}$ be the skew-adjacency matrix of $(H^\tau\ast G^\sigma)^o$.
The skew-spectrum of $(H^\tau\ast G^\sigma)^o$ is determined in the following theorem.
\begin{thm}\label{spectrum2}
Let $H^\tau$ be an oriented bipartite graph of order $m$ and let the
skew-eigenvalues of $H^\tau$ be the non-zero values $\pm \mu_1i$, $\pm
\mu_2i$, $\dots$, $\pm \mu_ti$ and $m-2t$ $0$'s. Let $G^\sigma$ be an oriented
graph of order $n$ and let the skew-eigenvalues of $G^\sigma$ be the non-zero
values $\pm \lambda_1i$, $\pm\lambda_2i$, $\dots$, $\pm \lambda_ri$ and $n-2r$ $0$'s.
Then the skew-eigenvalues of the oriented graph $(H^\tau\ast G^\sigma)^o$ are
$\pm\, i\sqrt{(u_j^2+1)(\lambda_k^2+1)-1}$ with multiplicities $2$, $j=1,\dots,t$, $k=1,\dots,r$,
$\pm\mu_ji$ with multiplicities $n-2r$, $j=1,\dots,t$, $\pm\lambda_ki$ with
multiplicities $m-2t$, $k=1,\dots,r$, and $0$ with multiplicities $(m-2t)(n-2r)$.
\end{thm}
\pf Suppose that $(X,Y)$ is the bipartition of the vertices of $H$ with $X=m_1$ and $Y=m_2$. Let
$S_1$ and $S_2$ be the skew-adjacency matrices of $H^\tau$ and $G^\sigma$, respectively, 
where $S_1=\left(\begin{array}{cc} 0& A\\ -A^T& 0\end{array}\right)$ and $A$ is an $m_1\times m_2$ 
matrix. Then the skew-adjacency matrix $\widehat{S}$ of $(H^\tau\ast G^\sigma)^o$ can be written as 
$\widehat{S}=\overline{S}+S$, where $\overline{S}$ and $S$ are the skew-adjacency matrices of 
$(H^\tau\square G^\sigma)^o$ and $(H^\tau\otimes G^\sigma)^o$, respectively.

With suitable labeling of the vertices of $H\ast G$, we can derive the following formulas.
\begin{equation}\label{equ}
\overline{S}=I'_{m_1+m_2}\otimes S_2+S_1\otimes I_n \hspace{10pt}\text{and}\hspace{10pt} S=S'_1\otimes S_2,
\end{equation}
where $I'_{m_1+m_2}=\left(\begin{array}{cc} I_{m_1}& 0\\ 0& -I_{m_2}\end{array}\right)$ and  
$S'_1=\left(\begin{array}{cc} 0& A\\ A^T& 0\end{array}\right)$. For the details of 
Equation (\ref{equ}), one can also see Theorem $3.1$ of \cite{CLL3} and Theorem \ref{spectrum} of this paper.

We then compute the singular values of $\widehat{S}$. Note that $\widehat{S}\,\widehat{S}^T=\overline{S}\,\,\overline{S}^T+SS^T+\overline{S}S^T+S\overline{S}^T$.
From Theorem $3.1$ of \cite{CLL3} or direct computation, we can derive that $\overline{S}\,\,\overline{S}^T=-(I_m\otimes S_2^2+S_1^2\otimes I_n)$. It is obvious that
$SS^T=(S'_1\otimes S_2)((S'_1)^T\otimes S_2^T)=-(S'_1)^2\otimes S_2^2=S_1^2\otimes S_2^2$.
Moreover,
\begin{equation*}
\begin{split}
\overline{S}S^T+S\overline{S}^T=&\left(I'_{m_1+m_2}\otimes S_2+S_1\otimes I_n\right)\left(S'_1\otimes (-S_2)\right)\\
&\hspace{5pt}+\left(S'_1\otimes S_2\right)\left(I'_{m_1+m_2}\otimes (-S_2)+(-S_1)\otimes I_n\right)\\
=&-\left[\left(S_1\otimes S_2^2+S_1S'_1\otimes S_2\right)+\left((-S_1)\otimes S_2^2+S'_1S_1\otimes S_2\right)\right]\\
=&0.
\end{split}
\end{equation*}

To sum up all computation, we obtain that
$$\widehat{S}\,\widehat{S}^T=-I_m\otimes S_2^2-S_1^2\otimes I_n+S_1^2\otimes S_2^2=(S_1^2-I_m)\otimes(S_2^2-I_n)-I_{mn}.$$

Therefore, the eigenvalues of $\widehat{S}\,\widehat{S}^T$ are $(\mu^2+1)(\lambda^2+1)-1$,
where $\mu i\in Sp(H^\tau)$ and $\lambda i\in Sp(G^\sigma)$. Then the skew-spectrum of
$(H^\tau\ast G^\sigma)^o$ immediately follows. The proof is complete. \qed

Similar to Theorem \ref{max-fami1}, we can construct a new family of oriented graphs with the
maximum skew energy by applying the above theorem.
\begin{thm}\label{max-fami2}
Let $H^\tau$ be an oriented $k$-regular bipartite graph of order $m$ with maximum
skew energy $m\sqrt{k}$. Let $G^\sigma$ be an oriented $l$-regular graph of order $n$
and maximum skew energy $n\sqrt{l}$. Then $(H^\tau\ast G^\sigma)^o$ is an oriented
$(k+l+kl)$-regular graph and has the maximum skew energy $\mathcal{E}_S((H^\tau\ast G^\sigma)^o)=mn\sqrt{k+l+kl}$.
\end{thm}

Comparing Theorems \ref{max-fami1}, \ref{max-fami2} obtained above with Theorem \ref{Cartesian}
(or see Theorem $3.2$ in \cite{CLL3}), we find that the oriented graphs constructed from these
theorems have the same order $mn$ but different regularities, which are $kl$, $k+l+kl$ and $k+l$, respectively.

\begin{figure}[h,t,b,p]
\begin{center}
\scalebox{1}[1]{\includegraphics{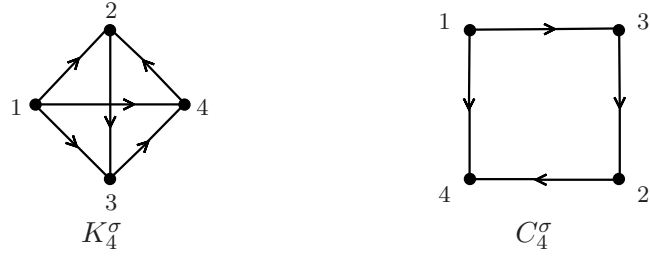}}
\end{center}
\caption{The orientations of $K_4$ and $C_4$ with the maximum skew energies}\label{Fig1}
\end{figure}

\begin{exam}
Let $H=C_4$, $G_0=K_4$, $G_1=H\square G_0$, \ldots, $G_r=H\square G_{r-1}$. Obviously, $G_r$ is a $(2r+3)$-regular graph of order $4^{r+1}$. From \cite{ABC,Tian}, we know that
$H$ has the orientation with the maximum skew energy $4\sqrt{2}$ and $G_0$ has the orientation with the maximum skew energy $4\sqrt{3}$, see Figure \ref{Fig1}. By Theorem \ref{Cartesian}, $G_r$ has the orientation with the maximum skew energy $4^{r+1}\sqrt{2r+3}$.
\end{exam}

\begin{exam}\label{exam1}
Let $H=C_4$, $G_0=K_4$, $G_1=H\otimes G_0$, \ldots, $G_r=H\otimes G_{r-1}$. It is obvious that $G_r$ is a 
$(3\cdot2^r)$-regular graph of order $4^{r+1}$. Then by Theorem \ref{max-fami1}, $G_r$ has the orientation 
with maximum skew energy $4^{r+1}\sqrt{3\cdot2^r}$.\\
\end{exam}

\begin{exam}\label{exam2}
Let $H=C_4$, $G_0=K_4$, $G_1=H\otimes G_0$, $G_2=G_1\otimes G_0$, \ldots, $G_r=G_{r-1}\otimes G_0$. 
Note that $H$, $G_1$, $G_2$, \dots, $G_{r-1}$ are all regular bipartite graphs and 
$G_r$ is a $(2\cdot 3^r)$-regular bipartite graph of order $4^{r+1}$. Then by Theorem \ref{max-fami1}, 
$G_r$ has the orientation with maximum skew energy $4^{r+1}\sqrt{2\cdot 3^r}$.
\end{exam}

\begin{exam}\label{exam3}
Let $H=C_4$, $G_0=K_4$, $G_1=H\ast G_0$, \ldots, $G_r=H\ast G_{r-1}$. Note that $G_r$ is a 
$(4\cdot 3^r-1)$-regular graph of order $4^{r+1}$. Then by Theorem \ref{max-fami2}, $G_r$ has 
the orientation with maximum skew energy $4^{r+1}\sqrt{4\cdot 3^r-1}$.
\end{exam}

From Examples \ref{exam1}, \ref{exam2} and \ref{exam3}, we can see that for some positive integers $k$, 
there exist oriented $k$-regular graphs with the maximum skew energy, which has order $n\leq k^2$. 
It is unknown that whether for any positive integer $k$, the oriented graph exists such that its 
order $n$ is less than $k^2$ and it has an orientation with the maximum skew energy.

\subsection{The orientation of $H[G]$}\label{lexi}

In this subsection, we consider the lexicographic product $H[G]$ of a bipartite graph $H$ and a graph $G$.
All definitions and notations are the same as above. We can see that the edge set $H[G]$ is the disjoint-union of the edge sets of $H\square G$ and $H\otimes K_n$, where $K_n$ is a complete graph of order $n$.

Let $H^\tau$ and $G^\sigma$ be oriented graphs of $H$ and $G$ with the skew-adjacency matrices $S_1$ and $S_2$, respectively. Let $K_n^\varsigma$ be an oriented graph of $K_n$ with the skew-adjacency matrix $S_3$. Then we can obtain two oriented graphs $(H^\tau\square G^\sigma)^o$ and $(H^\tau\otimes K_n^\varsigma )^o$. Thus it is natural to yield an orientation of $H[G]$, denoted by $H[G]^o$, such that the arc set of $H[G]^o$ is the disjoint-union of the arc sets  of $(H^\tau\square G^\sigma)^o$ and $(H^\tau\otimes K_n^\varsigma )$. Let $\widetilde{S}$ be the skew-adjacency matrix of $H[G]^o$. We can see that
$\widetilde{S}=\overline{S}+S'_1\otimes S_3$, where $\overline{S}=I'_{m_1+m_2}\otimes S_2+S_1\otimes I_n$ is the skew-adjacency matrix of $(H^\tau\square G^\sigma)^o$. Then
\begin{equation*}
\begin{split}
\widetilde{S}^T\widetilde{S}&=-\widetilde{S}^2=-\left(\overline{S}+S'_1\otimes S_3\right)^2\\
&=-\left[I_m\otimes S_2^2+S_1^2\otimes I_n+(S'_1)^2\otimes S_3^2+\overline{S}(S'_1\otimes S_3)+(S'_1\otimes S_3)\overline{S}\right]\\
&=-\left[I_m\otimes S_2^2+S_1^2\otimes I_n-S_1^2\otimes S_3^2+S_1\otimes (S_2S_3)-S_1\otimes (S_3S_2)\right]
\end{split}
\end{equation*}

Suppose that $H^\tau$ is an oriented $k$-regular bipartite graph of order $m$ with maximum skew energy $m\sqrt{k}$. Then $S_1^TS_1=kI_m$. Let $G^\sigma$ be an oriented $l$-regular graph of order $n$
and maximum skew energy $n\sqrt{l}$. Then $S_2^TS_2=lI_n$. It is obvious that $H[G]^o$ is $(kn+l)$-regular. Moreover, let $K_n^\varsigma$ be an oriented graph of $K_n$ with maximum skew energy $n\sqrt{n-1}$. Then $S_3^TS_3=(n-1)I_n$, that is, $S_3$ is a skew-symmetric Hardamard matrix \cite{MS} of order $n$. If another condition that $S_2S_3=S_3S_2$ holds, then
$$\widetilde{S}^T\widetilde{S}=-\left[I_m\otimes S_2^2+S_1^2\otimes I_n-S_1^2\otimes S_3^2\right]=(kn+l)I_{mn}.$$
By Lemma \ref{max}, $H[G]^o$ has the maximum skew energy $mn\sqrt{kn+l}$.

The following example illustrates that the oriented graph satisfying the above conditions indeed exists.
\begin{exam}
Let $H^\tau$ is an arbitrary oriented $k$-regular bipartite graph of order $m$ with the maximum skew energy $m\sqrt{k}$. Let $C_4^\sigma$ be the oriented graph of $C_4$ with maximum skew energy $4\sqrt{2}$ and the skew-adjacency matrix $S_2$, and $K_4^\varsigma$ be the oriented graph of $K_4$ with maximum skew energy $4\sqrt{3}$ and the skew-adjacency matrix $S_3$, see Figure \ref{Fig1}. It can be verified that $S_2S_3=S_3S_2$. It follows that $(H[G])^o$ is an oriented $(4k+2)$-regular graph of order $4m$ with maximum skew energy $4m\sqrt{4k+2}$.

There are many options for $H$, such as $P_2$, $C_4$, $K_{4,4}$, the hypercube $Q_d$ and so on, which forms a new family of oriented graphs with the maximum skew energy.
\end{exam}


\begin{thebibliography}{20}

\bibitem{ABC}C. Adiga, R. Balakrishnan, W. So, The shew energy of a digraph, Linear Algebra Appl. 432(2010), 1825--1835.

\bibitem{CLL}X. Chen, X. Li, H. Lian, $4$-Regular oriented graphs with optimum skew energy,
    arXiv: 1304.0847.

\bibitem{CLL3}X. Chen, X. Li, H. Lian, More on the skew-spectra of bipartite graphs and Cartesian products of graphs, arXiv: 1305.3414.

\bibitem{CH}D. Cui, Y. Hou, On the skew spectra of Cartesian products of graphs, Electron. J. Combin.
    20(2013), \#P19.

\bibitem{GX}S. Gong, G. Xu, $3$-Regular digraphs with optimum skew energy, Linear Algebra Appl. 436(2012), 465--471.

\bibitem{G}I. Gutman, The energy of a graph, Ber. Math. Statist. Sekt. Forschungsz. Graz, 103(1978), 1--22.

\bibitem{LL}X. Li, H. Lian, A survey on the skew energy of oriented graphs, arXiv: 1304.5707.

\bibitem{MS}F.J. MacWilliams, N.J.A. Sloane, The Theory of Error-Correcting Codes, North-Holland, New York, 1977.

\bibitem{Tian}G. Tian, On the skew energy of orientations of hypercubes, Linear Algebra Appl. 435(2011), 2140--2149.

\end{thebibliography}
\end{document}